\documentclass[reqno,preprint,12pt]{elsarticle}
\usepackage[active]{srcltx}
\usepackage{bbm}
\usepackage{tipa}
\usepackage{amsthm}
\usepackage{amssymb}
\usepackage{stackrel}
\usepackage{mathrsfs}
\usepackage{amsmath,amscd}
\usepackage{lscape}
\usepackage{longtable}
\usepackage{booktabs}
\usepackage{rotating}
\usepackage[active]{srcltx}
\usepackage{amsfonts}
\usepackage[centerlast]{caption}
\usepackage{latexsym}
\usepackage{amsmath,amsthm}
\usepackage{amssymb,amsmath}
\usepackage[all]{xy}
\usepackage{graphicx}
\usepackage{rotating}
\usepackage[dvips]{geometry}
 \newtheorem{dl}{Theorem}

\newtheorem{yl}[dl]{Lemma}
\newtheorem{dy}[dl]{Definition}

\newtheorem{remark}[dl]{Remark}

\numberwithin{equation}{section}
\newcommand{\bbtm}[4]{\bibitem{kn:#1}{#2,}~\emph{#3,}~{#4.}}
\newcommand{\cito}[1]{\cite{kn:#1}}
\newcommand{\citu}[2]{\cite[#2]{kn:#1}}
\newproof{pot1}{Proof of Theorem \ref{ma1}}
\newproof{pot2}{Proof of Theorem \ref{ma2}}
\def\qed{\hfill \rule{4pt}{7pt}}
\def\pf{\noindent {\sl Proof.~}}
\newcommand{\poq}[2]{(#1;q)_{#2}}

\usepackage[colorlinks=true,linkcolor=green,filecolor=brown,citecolor=green]{hyperref}
\usepackage{hyperref}
\usepackage{epsf}

\begin{document}
\title{ A new recurrence relation for the truncated very-well-poised $_6\psi_6$ series and Bailey's  summation formula}
\author{Jin Wang\fnref{fn3,fn4}}
\fntext[fn3]{Supported by NSF of Zhejiang Province (Grant~No.~LQ20A010004) and NSF of China (Grant~No. 12001492).}
\fntext[fn4]{E-mail address: jinwang@zjnu.edu.cn}
\address[P.R.China]{Department of Mathematics, Zhejiang Normal  University,
Jinhua 321004,~P.~R.~China}
\author{Xinrong Ma\fnref{fn1,fn2}}
\fntext[fn1]{Supported by NSFC grant No. 11971341.}
\fntext[fn2]{Corresponding author. E-mail address: xrma@suda.edu.cn.}
\address[P.R.China]{Department of Mathematics, Soochow University, Suzhou 215006, P.R.China}
\begin{abstract}
In this paper we introduce the truncated very-well-poised $_6\psi_6$ series and set up an explicit  recurrence relation for it by means of  the  classical Abel lemma on summation by parts. This new recurrence relation implies an   elementary  proof of  Bailey's
 well-known  $_6\psi_6$ summation formula.
\end{abstract}
\begin{keyword} Basic hypergeometric series; truncated; very-well-poised;          Bailey's  $_6\psi_6$ summation formula; Abel's lemma; Weierstrass'  theta  identity\\

{\sl AMS subject classification 2000}:  05A10; 33D15
\end{keyword}
\maketitle
\thispagestyle{plain}
\markboth{J. Wang and X. R. Ma}{A new recurrence relation for the truncated    VWP ${}_6\psi_6$ series and Bailey's summation formula}

\section{Introduction}

As is well known,
Bailey's bilateral very-well-poised (in short, VWP) ${}_6\psi_6$ summation formula is one of the deepest results in the theory of basic hypergeometric series,
which can be recorded as follows.
\begin{dl}[Bailey's VWP  ${}_6\psi_6$ summation formula: \citu{10}{(II.33)}]\label{baileypsi66} Let $a,b,c,d,e$  be
five nonzero complex parameters subject to $|a^2q/(bcde)|<1$. Then there holds
the summation formula
\begin{align}
&{}_6\psi_6\left(a;b,c,d,e
; q, \frac{a^2q}{bcde}\right)\label{bailey-a}\\
&=
\frac{(q,aq,q/a,aq/(be),aq/(ce),aq/(de),aq/(bc),aq/(bd),aq/(cd);q)_\infty}
{(aq/b,aq/c,aq/d,aq/e,q/b,q/c,q/d,q/e,a^{2}q/(bcde);q)_\infty}.\nonumber
\end{align}
\end{dl}
We would like to refer the reader to \cite{kn:andrews} by G. E. Andrews for some applications of the ${}_6\psi_6$ summation to
partitions and number theory.

To the best of our knowledge, finding simple and elementary proof of Bailey's  VWP $_6\psi_6$ summation formula is still one of  attractive problems in basic hypergeometric series.  As a supporting evidence, we would like to readdress the comment of R. Askey in his paper \cite[p. 575]{kn:askey0} ``\ldots \textit{However, it is annoying that a sum that is this important has not been obtained from a more elementary
special case.}”    Up to now,  many different proofs of Bailey's  VWP $_6\psi_6$ summation formula have been found,  among are the method of  integral and functional equations by R. Askey \cite{kn:askey0}, the method of  Liouville’s analytic continuation   by R. Askey and M. E. H. Ismail \cite{kn:askey}, the $q$-difference method with series expansion by G. E.  Andrews \cite{kn:andrews}, the method of $q$-Gosper algorithm by V. Y. B. Chen, W. Y. C. Chen, and N. S. S. Gu \cite{kn:chengu},  the difference method together Abel's lemma by W. Chu \cite{kn:chu0}, the method of Cauchy's residue by  F. Jouhet and M. Schlosser \cite{kn:schlosser2,kn:ttt}, and  elementary manipulations of series by M. Schlosser alone \cite{kn:ttt},  L. J. Slater and A. Lakin \cite{kn:slater}.  It should be  mentioned that   in his paper \cite{kn:bailey1936}, Bailey described  how to deduce Bailey's VWP ${}_6\psi_6$ summation formula \citu{bailey1936}{(4.7)} and Weierstrass' theta  identity \citu{bailey1936}{(5.2)} from some three-term relations for VWP ${}_8\phi_7$ series, provided that Rogers'  VWP ${}_6\phi_5$ summation formula is given. However, he did not give any direct connection between these summation formulas. In our paper \cite{kn:wangma}, we have revealed
 certain relation among Bailey's VWP ${}_6\psi_6$ and Rogers'  ${}_6\phi_5$ summation formulas, as well as Weierstrass' theta identity.

For purpose of comparison,  we  especially point out that  it is just Abel's lemma on summation by parts with which W. Chu rediscovered  in a series of papers  such as \cite{kn:chu0,kn:chu2} many important results for basic hypergeometric series. Among these, there are  the $q$-binomial theorem, $q$-Gauss  theorem, Ramanujan's bilateral ${}_{1} \psi_{1}$, the $q$-Pfaff-Saalschütz sum, and Jackson's VWP ${}_{8} \phi_{7}$ sum, and Bailey's  VWP $_6\psi_6$ summation formula. In our view, both Chu's proof in \cite{kn:chu0} and Chen-Chen-Gu's proof in \cite{kn:chengu} use more than four recurrence relations and require the Jacobi triple product identity.

In this short article, as a possibly desired way by R. Askey, we will introduce
\begin{dy}
For any integer $N\geq 0$ and  nonzero complex parameters $A,B,C,D,E$, define the truncated very-well-poised $_6\psi_6$ series $\mathbf{S}_N(A,B,C,D,E)$ to be the following finite sum
 \begin{align}
\sum_{n=-N}^{N} \frac{\nabla\left(BDEq^{2n+1}/A\right)}{\nabla\left(BDEq/A\right)}\frac{\poq{Bq,  Dq, Eq,BCDEq^2/A^2}{n}}{
\poq{DEq/A,BEq/A, BDq/A, A/C}{n}}\bigg(\frac{1}{Cq^2}\bigg)^n,\label{SNsequence}
\end{align}
where the notation $\nabla(x):=1-x$.
\end{dy}
In our working below we often write $\mathbf{S}_N(A;C)$ for  $\mathbf{S}_N(A,B,C,D,E)$ by suppressing the dependence of the various summations on the complex parameters $B,D,E$ for easy of notation.  As one of the main results, we will  present  a somewhat more ``unexpected"  recurrence relation of \eqref{SNsequence} underlying  Bailey's VWP $_6\psi_6$ summation formula. In other word,  Bailey's VWP $_6\psi_6$ summation formula is just a limitation of this new recurrence relation. Our argument,  apart from  Abel's lemma on summation by parts, only depends on the following self-evident  identity.
\begin{yl}\label{theta}
For any  complex parameters $b,c,x,z$ with $bcz\neq 0$,  we have
\begin{align}
\nabla\left(cx,\frac{x}{c},bz,
\frac{z}{b}\right)-\nabla \left(bx,\frac{x}{b},cz,\frac{z}{c}\right)
=\frac{z}{c}\nabla \left(bc,\frac{c}{b},
xz,\frac{x}{z}\right).\label{trivalweierstrass-new-1}
\end{align}
Hereafter, for brevity,  we employ  the notation
\begin{align}
\nabla(x_1,x_2,\ldots,x_n):=\prod_{i=1}^n\nabla(x_i).\label{nabladef}
\end{align}
\end{yl}

Some remarks on notation are necessary.  Throughout this paper, we will adopt the standard notation and terminology for basic
hypergeometric series (or $q$-series) from the book  \cite{kn:10} (Gasper and Rahman, 2004). For instance, the $q$-shifted factorial with $0<|q|<1$ is defined by
\begin{align}
(a; q)_{n}:=\left\{\begin{array}{ll}
(1-a)(1-a q) \cdots\left(1-a q^{n-1}\right), & n=1,2, \ldots, \\
1, & n=0, \\
{\left(\left(1-a q^{-1}\right)\left(1-a q^{-2}\right) \cdots\left(1-a q^{n}\right)\right)^{-1},} & n=-1,-2, \ldots
\end{array}\right.
\end{align}
Its multi-parameter form is compactly abbreviated to
\[(x_1,x_2,\cdots,x_m;q)_n:\:=\:(x_1;q)_n(x_2;q)_n\cdots(x_m;q)_n.\]

The basic and bilateral hypergeometric series with the
base $q$ and the argument $z$ are defined, respectively, by
\begin{align}{}_{r}\phi _{r-1}\left[\begin{matrix}a_{1},\dots ,a_{r}
\\ b_{1},\dots ,b_{r-1}\end{matrix}
; q, z\right]&:\:=\:\sum _{n=0} ^{\infty }\frac{\poq {a_{1},\cdots
,a_{r}}{n}}{\poq
{q,b_{1},\cdots,b_{r-1}}{n}}z^{n},\label{phino}\\
{}_{r}\psi _{r}\left[\begin{matrix}a_{1},\dots ,a_{r}\\
b_{1},\dots ,b_{r}\end{matrix} ; q, z\right]&:\:=\:\sum _{n=-\infty}
^{\infty }\frac{\poq {a_{1},\cdots ,a_{r}}{n}}{\poq
{b_{1},\cdots,b_{r}}{n}}z^{n}.\label{yes}
\end{align}
In particular, the compact notation
$\,_{r}\psi_{r}(a;a_3,\cdots,a_{r};q,z)$ denotes the special case
 of  ${}_{r}\psi _{r}$ series above, called  \emph{very-well-poised} 
 (VWP), in which all parameters
satisfy the relations
\begin{align}
b_1a_1=b_2a_2=\cdots=b_{r}a_{r}=aq;~ a_1=q\sqrt{a},a_2=-q\sqrt{a}\label{conditions}
\end{align}
and
$\,_{r}\phi_{r-1}(a;a_4,\cdots,a_{r};q,z)$ if there exists certain  parameter $a_i$ (say $a_3$) $=a$ in \eqref{conditions}.

\section{A recurrence relation for $\mathbf{S}_N(A;C)$}

Let us begin with Abel's lemma on summation by parts.
 \begin{yl}[Abel's lemma]\label{ablespecial} For any two sequences $\{U_n, V_n|n=0,\pm 1,\pm 2,\cdots\}$ and integers $M,N\geq 0$, it always holds
\begin{align}\sum_{n=-M}^{N} V_n(U_n-U_{n+1})
=V_{-M}U_{-M}-V_NU_{N+1}+\sum_{n=-M+1}^{N} U_n(V_n-V_{n-1}).\label{abel}
\end{align}
\end{yl}

As one of our main results,  the following recurrence relation for $\mathbf{S}_N(A;C)$ may serve as an essential characteristic for  Bailey's VWP ${}_6\psi_6$ summation formula. It is a direct application of  Lemmas \ref{theta} and  \ref{ablespecial}.
\begin{dl}\label{bilateralseries}Let $\mathbf{S}_N(A;C)$ be defined by \eqref{SNsequence}.
Then  $\mathbf{S}_N(A;C)$ satisfies the following recurrence relation
\begin{align}
\mathbf{S}_{N+1}(A;C)
&=\frac{K_N(A;B,C,D,E)}{(Cq^3)^{N}}\label{bailey-rec}\\
&+\frac{A^2q}{BDE}\frac{\nabla(BDE/A,Cq^3,BD/A,BE/A, DE/A)}{\nabla\left(BDEq/A,BCDEq/A^2,Aq/B,Aq/D,Aq/E\right)}\mathbf{S}_{N}(Aq;Cq),\nonumber
\end{align}
where
 \begin{align}
&K_N(A;B,C,D,E):=\frac{q^{N-2}~\nabla\left(BDEq^{2N+3}/A\right)\poq{Bq,  Dq, Eq,BCDEq^2/A^2}{N+1}}{C~\nabla\left(BDEq/A\right)
\poq{A/C,BDq/A,BEq/A, DEq/A}{N+1}}\label{KN12}\\
&+\frac{BDE~\nabla\left(
q^{N-1}/A\right)\poq{A/BD,A/BE, A/DE,C/A}{N+2}}
  {C~\nabla\left(C/A,Aq/B,Aq/D,Aq/E,BDEq/A\right)\poq{1/B,  1/D, 1/E,A^2/(BCDEq)}{N+1}}\nonumber\\  &-\frac{A^2q}{BDE}\frac{\nabla\left(Bq,Dq,Eq,BDEq^{N-1}/A^{2}\right)\poq{BCDEq^2/A^2,Bq^2,  Dq^2, Eq^2}{N}}
{\nabla\left(Aq/B,Aq/D,Aq/E,BDEq/A\right)\poq{A/C,BDq/A,
BEq/A, DEq/A}{N}}\nonumber.
\end{align}
\end{dl}
\pf ~~To establish \eqref{bailey-rec}, we  start with two sequences
 \begin{subequations}\label{UVseqs}
\begin{align}
U_{n}&:=\frac{\poq{Bq,  Dq, Eq, BDE/A^{2}q}{n}}{\poq{BD/A,
BE/A, DE/A, Aq^2}{n}},
 \\
V_{n}&:=\frac{\poq{Aq^2,BCDEq/A^2}{n+1}}{
\poq{A/Cq,BDE/A^2q^2}{n+1}}\left(\frac{1}{Cq^3}\right)^{n}.
\end{align}
\end{subequations} In view of Lemma \ref{ablespecial},
we  need only to calculate  the differences
$U_n-U_{n+1}$ and $V_n-V_{n-1}.$
To do this, we start with  \eqref{trivalweierstrass-new-1}
 and make the parameter replacement
\begin{align}(b,c,x,z)\to\left(\frac{A q}{\sqrt{B} \sqrt{D}}, \frac{\sqrt{B}}{\sqrt{D}},\sqrt{B}  \sqrt{D}q,\frac{\sqrt{B} \sqrt{D} E}{A}\right).\label{ggggg}
\end{align}
As a result, it follows immediately
\begin{align*}
\nabla\left(cx,\frac{x}{c},bz,
\frac{z}{b}\right)&=\nabla\left(Bq,Dq,Eq,\frac{BDE}{A^2q}\right),\\
\nabla \left(bx,\frac{x}{b},cz,\frac{z}{c}\right)&=\nabla \left(Aq^2,\frac{BD}{A},\frac{BE}{A},\frac{DE}{A}\right),\\
\nabla \left(bc,\frac{c}{b},
xz,\frac{x}{z}\right)&=\nabla \left(\frac{A q}{D},\frac{B}{A q},\frac{BDEq}{A},\frac{A q}{E}\right),
\end{align*}
specializing \eqref{trivalweierstrass-new-1} to the form
\begin{align}
\nabla\left(Bq,Dq,Eq,\frac{BDE}{A^2q}\right)&-\nabla \left(A q^2,\frac{BD}{A},\frac{BE}{A},\frac{DE}{A}\right)\nonumber\\
&=\frac{DE}{A}\nabla \left(\frac{A q}{D},\frac{B}{A q},\frac{BDEq}{A},\frac{A q}{E}\right).\label{weierstrass-newnew-new-1}
\end{align}
In such case, it is easy to find
\begin{align*}
U_{n}-U_{n+1}&=U_{n}\times\bigg(1-\frac{U_{n+1}}{U_{n}}\bigg)\\
&=U_{n}\times\bigg(1-
\frac{\nabla\left(Bq^{1+n},Dq^{1+n},
Eq^{1+n},\frac{BDEq^{n}}{A^2q}\right)}
{\nabla \left(A q^{n+2},\frac{BD}{A}q^n,\frac{BE}{A}q^n,\frac{DE}{A}q^n\right)}\bigg)\\
&=U_{n}\times \bigg(1-
\frac{\nabla\left(Bq,Dq,
Eq,\frac{BDE}{A^2q}\right)}
{\nabla \left(Aq^2,\frac{BD}{A},\frac{BE}{A},\frac{DE}{A}\right)}\bigg)\bigg|_\sigma.
\end{align*}
Herein and in what follows, the notation $\sigma$
 denotes the parameter replacement
 \[\sigma: (A,B,D,E)\to(Aq^{n},Bq^n,Dq^n,Eq^n)\]
and $F\big|_\sigma$ means applying  $\sigma$ to the function $F$. Applying \eqref{weierstrass-newnew-new-1} to the right-hand side of the last identity  gives rise to
\begin{align}
U_{n}-U_{n+1}&=U_{n}\times \bigg(\frac{-DE}{A}
\frac{\nabla \left(\frac{A q}{D},\frac{B}{A q},\frac{BDEq}{A},\frac{A q}{E}\right)}
{\nabla \left(Aq^2,\frac{BD}{A},\frac{BE}{A},\frac{DE}{A}\right)}\bigg)\bigg|_\sigma\nonumber\\
&=U_{n}\times \frac{-DEq^n}{A}
\frac{\nabla\left(B/Aq,Aq/D,Aq/E,BDEq^{2n+1}/A\right)}
{\nabla \left(Aq^{n+2},\frac{BD}{A}q^n,\frac{BE}{A}q^n,\frac{DE}{A}q^n\right)}\nonumber\\
&=-\frac{DEq^n}{A}\nabla\left(B/Aq,Aq/D,Aq/E\right)\label{Udiffence}\\
&\quad\times\frac{\nabla \left(BDEq^{2n+1}/A\right)\poq{Bq,  Dq, Eq, BDE/A^{2}q}{n}}{\poq{BD/A,
BE/A, DE/A, Aq^2}{n+1}}.\nonumber
\end{align}
On the same lines, we can compute the difference  $V_n-V_{n-1}$. We still
 appeal to the basic identity \eqref{trivalweierstrass-new-1}
 and make, instead of \eqref{ggggg}, the following parameter replacement
$$(b,c,x,z)\to\left(\frac{A^{3/2}q^{1/2} }{\sqrt{B} \sqrt{C} \sqrt{D} \sqrt{E}},\frac{\sqrt{B}  \sqrt{D} \sqrt{E} q^{-3/2}}{\sqrt{A} \sqrt{C}},0,\frac{\sqrt{B} \sqrt{C} \sqrt{D} \sqrt{E}q^{3/2} }{\sqrt{A}}\right).
$$
In the  sequel, we have
\begin{align}\nabla\left(A q^2,\frac{BCDE q}{A^2}\right)-Cq^3\nabla \left(\frac{A}{Cq},\frac{BDE}{A^2q^2}\right)=\nabla \left(\frac{B D E}{A},Cq^3\right).\label{Vdiffence-two}
\end{align}
Now we proceed to calculate the difference
\begin{align*}
V_{n}-V_{n-1}&=-V_{n-1}\times\bigg(1-\frac{V_{n}}{V_{n-1}}\bigg)\\
&=-V_{n-1}\times\bigg(1-\frac{\nabla(Aq^{n+2},BCDEq^{n+1}/A^2)}{
\nabla(Aq^{n}/Cq,BDEq^{n}/A^2q^2)}\frac{1}{Cq^3}\bigg)\\
&=-V_{n-1}\times\bigg(1-\frac{\nabla(Aq^2,BCDEq/A^2)}{
\nabla(A/Cq,BDE/A^2q^2)}\frac{1}{Cq^3}\bigg)\bigg|_\tau,
\end{align*}
where $\tau$ stands for the parameter replacement
\[\tau: (A,B,C,D,E)\to(Aq^{n},Bq^{n},C,Dq^{n},Eq^{n}).\]
At this stage, by applying \eqref{Vdiffence-two} to the last identity, we arrive at
\begin{align}
V_{n}-V_{n-1}&=V_{n-1}\times\bigg(\frac{\nabla(BDE/A,Cq^3)}{
\nabla(A/Cq,BDE/A^2q^2)}\frac{1}{Cq^3}\bigg)\bigg|_\tau\nonumber\\
&=\frac{\poq{Aq^2,BCDEq/A^2}{n}}{
\poq{A/Cq,BDE/A^2q^2}{n}}\left(\frac{1}{Cq^3}\right)^{n-1}\times \frac{\nabla(BDEq^{2n}/A,Cq^3)}{
\nabla(Aq^{n}/Cq,BDEq^{n}/A^2q^2)}\frac{1}{Cq^3}\nonumber\\
&=\frac{\poq{Aq^2,BCDEq/A^2}{n}}{
\poq{A/Cq,BDE/A^2q^2}{n+1}}\times \frac{\nabla(BDEq^{2n}/A,Cq^3)}{(Cq^3)^n}. \label{Vdiffence}
\end{align}
This identity together with \eqref{Udiffence} specifies  \eqref{abel} to the form
\begin{align*}
&-\frac{DE}{A}\sum_{n=-M}^{N}\frac{\poq{Aq^2,BCDEq/A^2}{n+1}}{
\poq{A/Cq,BDE/A^2q^2}{n+1}}
\frac{\nabla\left(B/Aq,Aq/D,Aq/E,BDEq^{2n+1}/A\right)}{(Cq^2)^{n}}\\
&\qquad\qquad\times\frac{\poq{Bq,  Dq, Eq, BDE/A^{2}q}{n}}{\poq{BD/A,
BE/A, DE/A, Aq^2}{n+1}}\\
&=V_{-M}U_{-M}-V_NU_{N+1}+\sum_{n=-M+1}^{N} \frac{\poq{Bq,  Dq, Eq, BDE/A^{2}q}{n}}{\poq{BD/A,
BE/A, DE/A, Aq^2}{n}}\\
&\qquad\qquad\times\frac{\poq{Aq^2,BCDEq/A^2}{n}}{
\poq{A/Cq,BDE/A^2q^2}{n+1}}\frac{\nabla(BDEq^{2n}/A,Cq^3)}{(Cq^3)^n}.
\end{align*}
Simplifying the last identity by the relation  $\nabla(x)=-x\nabla(1/x),$ we have
\begin{align*}
&\frac{BDE}{A^2q}\frac{\nabla\left(Aq/B,Aq/D,Aq/E\right)}{\nabla(B,D,E,BDE/A^2q^2)}\\
&\qquad\times\sum_{n=-M}^{N}\frac{\nabla \left(BDEq^{2n+1}/A\right)}{(Cq^2)^n}\frac{\poq{B,  D, E,BCDEq/A^2}{n+1}}{
\poq{A/Cq,BD/A,BE/A, DE/A}{n+1}}\\
&=V_{-M}U_{-M}-V_NU_{N+1}+\frac{\nabla\left(Cq^3\right)}{\nabla(A/Cq,BDE/A^2q^2)}\\
&\qquad\times\sum_{n=-M+1}^{N}\frac{\nabla(BDEq^{2n}/A)}{(Cq^3)^n}\frac{\poq{BCDEq/A^2,Bq,  Dq, Eq}{n}}{\poq{A/C,BD/A,BE/A, DE/A}{n}}.
\end{align*}
By dividing both sides by
$$\frac{BDE}{A^2q}\frac{\nabla\left(Aq/B,Aq/D,Aq/E\right)}
{\nabla(B,D,E,BDE/A^2q^2)},$$
we obtain\begin{align}
&\sum_{n=-M}^{N} \frac{\nabla\left(BDEq^{2n+1}/A\right)}{(Cq^2)^n}\frac{\poq{B,  D, E,BCDEq/A^2}{n+1}}{
\poq{A/Cq,BD/A,BE/A, DE/A}{n+1}}\nonumber\\
=&(V_{-M}U_{-M}-V_NU_{N+1})\frac{A^2q}{BDE}\frac{\nabla(B,D,E,BDE/A^2q^2)}{\nabla\left(Aq/B,Aq/D,Aq/E\right)}\nonumber\\
&+\frac{A^2q}{BDE}\frac{\nabla(Cq^3,B,D,E)}
{\nabla\left(A/Cq,Aq/B,Aq/D,Aq/E\right)}\label{Important-happy}\\
&\times\sum_{n=-M+1}^{N}\frac{\nabla(BDEq^{2n}/A)}{(Cq^3)^n}\frac{\poq{BCDEq/A^2,Bq,  Dq, Eq}{n}}{\poq{A/C,BD/A,BE/A, DE/A}{n}}.\nonumber
\end{align}
Furthermore, on
 multiplying  both sides of \eqref{Important-happy} with
$$\frac{\nabla\left(A/Cq,BD/A,BE/A,DE/A\right)}{\nabla(B,D,E,BCDEq/A^2,BDEq/A)}
,$$
we obtain
\begin{align}
&\sum_{n=-M}^{N}
\frac{\nabla\left(BDEq^{2n+1}/A\right)}{\nabla\left(BDEq/A\right)}
\frac{\poq{Bq,  Dq, Eq,BCDEq^2/A^2}{n}}{
\poq{A/C,BDq/A,BEq/A, DEq/A}{n}}\bigg(\frac{1}{Cq^2}\bigg)^n\nonumber\\
=&(V_{-M}U_{-M}-V_NU_{N+1})\frac{A^2q}{BDE}
\frac{\nabla\left(A/Cq,BD/A,BE/A,DE/A,BDE/A^2q^2\right)}
{\nabla\left(Aq/B,Aq/D,Aq/E,BCDEq/A^2,BDEq/A\right)}\nonumber\\
&+\frac{A^2q}{BDE}\frac{\nabla(Cq^3,BD/A,BE/A,DE/A,BDE/A)}
{\nabla\left(Aq/B,Aq/D,Aq/E, BCDEq/A^2,BDEq/A\right)}\label{Important-happyhappy}\\
&\times\sum_{n=-M+1}^{N}\frac{\nabla(BDEq^{2n}/A)}{\nabla(BDE/A)}\frac{\poq{BCDEq/A^2,Bq,  Dq, Eq}{n}}{\poq{A/C,BD/A,BE/A, DE/A}{n}}\bigg(\frac{1}{Cq^3}\bigg)^n.\nonumber
\end{align}
The final step is to  calculate two terms $V_{-M}U_{-M}$ and $V_NU_{N+1}$. For this, we easily find
\begin{align*}
U_{N+1}V_{N}&=\frac{\poq{BCDEq/A^2,Bq,  Dq, Eq, BDE/A^{2}q}{N+1}}{\poq{A/Cq,BDE/A^2q^2,BD/A,
BE/A, DE/A}{N+1}}\left(\frac{1}{Cq^3}\right)^{N}
\end{align*}
while,  according to the basic relation (see \cite[(I.11)]{kn:10})
\begin{align}
\frac{(x ; q)_{-m}}{(y ; q)_{-m}}=\frac{(q / y ; q)_{m}}{(q / x ; q)_{m}}\left(\frac{y}{x}\right)^{m},\label{nongetive}
\end{align}
we easily check
\begin{multline*}
V_{-M}U_{-M}\\=\frac{\poq{Bq,  Dq, Eq, BDE/A^{2}q}{-M}}{\poq{BD/A,
BE/A, DE/A, Aq^2}{-M}}
 \times\frac{\poq{Aq^2,BCDEq/A^2}{-M+1}}{
\poq{A/Cq,BDE/A^2q^2}{-M+1}}\left(\frac{1}{Cq^3}\right)^{-M}\\
=\frac{\poq{Aq/BD,Aq/BE, Aq/DE, 1/Aq}{M}}{\poq{1/B,  1/D, 1/E, A^{2}q^2/BDE}{M}}
 \times
 \frac{
\poq{Cq^2/A,A^2q^3/BDE}{M-1}}{\poq{1/Aq,A^2/BCDE}{M-1}}\left(\frac{1}{Cq^3}\right)^{M-2}.
\end{multline*}
For our purpose, here we  need only to consider the case $M=N+1$. As such,  \eqref{Important-happyhappy} can be recast into the form
\begin{align}
&\sum_{n=-N-1}^{N+1} \frac{\nabla\left(BDEq^{2n+1}/A\right)}{\nabla\left(BDEq/A\right)}\frac{\poq{Bq,  Dq, Eq,BCDEq^2/A^2}{n}}{
\poq{A/C,BDq/A,BEq/A, DEq/A}{n}}\bigg(\frac{1}{Cq^2}\bigg)^n\nonumber\\
&-\frac{A^2q}{BDE}\frac{\nabla(Cq^3,BD/A,BE/A,DE/A,BDE/A)}
{\nabla\left(Aq/B,Aq/D,Aq/E, BCDEq/A^2,BDEq/A\right)}\nonumber\\
&\times\sum_{n=-N}^{N}\frac{\nabla(BDEq^{2n}/A)}{\nabla(BDE/A)}\frac{\poq{BCDEq/A^2,Bq,  Dq, Eq}{n}}{\poq{A/C,BD/A,BE/A, DE/A}{n}}\bigg(\frac{1}{Cq^3}\bigg)^n\nonumber\\
&=K^{(1)}_N(A;B,C,D,E)(Cq^3)^{-N}-K^{(2)}_N(A;B,C,D,E)(Cq^3)^{-N}\nonumber\\
&\qquad+K^{(3)}_N(A;B,C,D,E)(Cq^2)^{-N-1},\label{yeysyesyes}
\end{align}
where for $i=1,2,3,$~$K_N^{(i)}(A;B,C,D,E)$ are defined, respectively,  by
\begin{multline*}
K_N^{(1)}(A;B,C,D,E)\\
:=V_{-N-1}U_{-N-1}\frac{A^2q}{BDE}
\frac{\nabla\left(A/Cq,BD/A,BE/A,DE/A,BDE/A^2q^2\right)}
{\nabla\left(Aq/B,Aq/D,Aq/E,BCDEq/A^2,BDEq/A\right)}(Cq^3)^N;\end{multline*}
\begin{multline*}
K_N^{(2)}(A;B,C,D,E)\\:=V_{N}U_{N+1}\frac{A^2q}{BDE}
\frac{\nabla\left(A/Cq,BD/A,BE/A,DE/A,BDE/A^2q^2\right)}
{\nabla\left(Aq/B,Aq/D,Aq/E,BCDEq/A^2,BDEq/A\right)}(Cq^3)^N;
\end{multline*}
\begin{multline*}
K^{(3)}_N(A;B,C,D,E)\\
:=\frac{\nabla\left(BDEq^{2N+3}/A\right)}{\nabla\left(BDEq/A\right)}\frac{\poq{Bq,  Dq, Eq,BCDEq^2/A^2}{N+1}}{
\poq{A/C,BDq/A,BEq/A, DEq/A}{N+1}}.
\end{multline*}
For more clarity,  let us write $K_N(A;B,C,D,E)$ for the sum
\begin{align*}
K^{(1)}_N(A;B,C,D,E)-K^{(2)}_N(A;B,C,D,E)+K^{(3)}_N(A;B,C,D,E)q^{N-2}/C,
\end{align*}
which, after some routine computations, is the same as given by \eqref{KN12}.
In conclusion,  we are  able  to reformulate \eqref{yeysyesyes} in terms of the truncated VWP series $\mathbf{S}_N(A;C)$ and $K_N(A;B,C,D,E)$ as follows:
\begin{multline*}
\mathbf{S}_{N+1}(A;C)
=\frac{A^2q}{BDE}\frac{\nabla(BDE/A,Cq^3,BD/A,BE/A, DE/A)}{\nabla\left(BDEq/A,BCDEq/A^2,Aq/B,Aq/D,Aq/E\right)}\mathbf{S}_{N}(Aq;Cq)\\
+K_N(A;B,C,D,E)(Cq^3)^{-N}.
\end{multline*}
This gives the complete proof of Theorem \ref{bilateralseries}.
\qed

\section{A new proof of Bailey's VWP ${}_6\psi_6$ summation formula}
Having established Theorem  \ref{bilateralseries}, we now turn to show that  Bailey's VWP ${}_6\psi_6$ summation formula can be derived from the limit of the truncated VWP series $\mathbf{S}_N(A;C)$ as $N\to \infty$.  Actually, from Theorem  \ref{bilateralseries},  we may derive  without any difficulty that
\begin{yl} \label{MMMMM-YL}For  any
five nonzero complex parameters $B,C,D,E,X$  subject to $|C/q^3|<1$, define
\begin{align}
\mathbf{T}(X;C):={}_6\psi_6(BCDEX^2;BCDEXq,BXq,  DXq, EXq;q,C/q^3).
\end{align}
Then it holds
\begin{align}
\mathbf{T}(X;C)=\frac{\nabla\left(1/(BCDEXq),BCDEX^2q,BC/q,CD/q,CE/q\right)}
{\nabla(1/(BCDEX^2),C/q^3,BCDX,BCEX, CDEX)}\mathbf{T}(X;Cq).\label{MMMMM}
\end{align}
In particular,
\begin{align}
\mathbf{T}(q;C)=\frac{\poq{BCDEq^3,BC/q,CD/q,CE/q}{\infty}}
{\poq{C/q^3,BCDq,BCEq, CDEq}{\infty}}.\label{watson}
\end{align}
\end{yl}
\pf  Starting from Theorem  \ref{bilateralseries} with the tentative assumption $|Cq^2|>1$, we can take the limit of \eqref{bailey-rec} as  $N\to +\infty$. It is easy to check
$$
\lim_{N\to +\infty}\frac{K_N(A;B,C,D,E)}{(Cq^3)^N}=0,
$$
which  results from the fact
\begin{align*}
&\lim_{N\to +\infty}K_N(A;B,C,D,E)\\
=&\frac{BDE~\poq{A/BD,A/BE, A/DE,Cq/A}{\infty}}
  {C~\nabla\left(Aq/B,Aq/D,Aq/E,BDEq/A\right)\poq{1/B,  1/D, 1/E,A^2/(BCDEq)}{\infty}}\nonumber\\  &-\frac{A^2q}{BDE}\frac{\poq{BCDEq^2/A^2,Bq,Dq,Eq}{\infty}}
{\nabla\left(Aq/B,Aq/D,Aq/E,BDEq/A\right)\poq{A/C,BDq/A,
BEq/A, DEq/A}{\infty}}\\
&=\frac{BDE~}
  {C\nabla\left(Aq/B,Aq/D,Aq/E,BDEq/A\right)}\\
  &\times\frac{\theta(A/BD,A/BE, A/DE,A/C;q)-\frac{A^2Cq}{(BDE)^2}\theta(1/B,  1/D, 1/E,A^2/(BCDEq);q)}
  {\poq{1/B,  1/D, 1/E,A^2/(BCDEq),A/C,BDq/A,BEq/A, DEq/A}{\infty}}.
\end{align*}
All together, we conclude that
\begin{align*}
\mathbf{S}_{\infty}(A;C)=\frac{A^2q}{BDE}\frac{\nabla(BDE/A,Cq^3,BD/A,BE/A, DE/A)}{\nabla\left(BDEq/A,BCDEq/A^2,Aq/B,Aq/D,Aq/E\right)}~ \mathbf{S}_{\infty}(Aq;Cq),
\end{align*}
which, written out in full form, amounts to
\begin{align}
&\sum_{n=-\infty}^{\infty} \frac{\nabla\left(BDEq^{2n+1}/A\right)}{\nabla\left(BDEq/A\right)}
\frac{\poq{Bq,  Dq, Eq,BCDEq^2/A^2}{n}}{
\poq{A/C,BDq/A,BEq/A, DEq/A}{n}}\bigg(\frac{1}{Cq^2}\bigg)^n\nonumber\\
=&\frac{A^2q}{BDE}\frac{\nabla(BDE/A,Cq^3,BD/A,BE/A, DE/A)}{\nabla\left(BDEq/A,BCDEq/A^2,Aq/B,Aq/D,Aq/E\right)}\label{important}\\
&\times\sum_{n=-\infty}^{\infty} \frac{\nabla(BDEq^{2n}/A)}{\nabla(BDE/A)}\frac{\poq{BCDEq/A^2,Bq,  Dq, Eq}{n}}{\poq{A/C,BD/A,BE/A, DE/A}{n}}\bigg(\frac{1}{Cq^3}\bigg)^n.\nonumber
\end{align}
Next,  we further make the simultaneous substitution
\[(A,B,D,E)\to(CX,BX,DX,EX)\]
in \eqref{important} and then  replace $C$ with $1/C$. Consequently, we have
\begin{align}
&\sum_{n=-\infty}^{\infty} \frac{\nabla\left(BCDEX^2q^{2n+1}\right)}{\nabla\left(BCDEX^2q\right)}\frac{\poq{BCDEXq^2,BXq,  DXq, EXq}{n}}{
\poq{X,BCDXq,BCEXq, CDEXq}{n}}\bigg(\frac{C}{q^2}\bigg)^n\nonumber\\
=&\frac{q}{C^2BDEX}\frac{\nabla(BCDEX^2,q^3/C,BCDX,BCEX, CDEX)}{\nabla\left(BCDEX^2q,BCDEXq,q/BC,q/CD,q/CE\right)}\label{425}\\
&\times\sum_{n=-\infty}^{\infty} \frac{\nabla(BCDEX^2q^{2n})}{\nabla(BCDEX^2)}\frac{\poq{BCDEXq,BXq,  DXq, EXq}{n}}{\poq{X,BCDX,BCEX, CDEX}{n}}\bigg(\frac{C}{q^3}\bigg)^n.\nonumber
\end{align}
It is of importance to realize that the infinite sum on the far right-hand side of \eqref{425} is just $\mathbf{T}(X;C)$ while the left-hand side of \eqref{425} is nothing but $\mathbf{T}(X;Cq)$. Then we obtain the following recursive relation
\begin{align*}
\mathbf{T}(X;C)
=&\frac{X}{q}\frac{\nabla\left(BCDEX^2q,BCDEXq,BC/q,CD/q,CE/q\right)}
{\nabla(BCDEX^2,C/q^3,BCDX,BCEX, CDEX)}\mathbf{T}(X;Cq)
\\
=&\frac{\nabla\left(1/(BCDEXq),BCDEX^2q,BC/q,CD/q,CE/q\right)}
{\nabla(1/(BCDEX^2),C/q^3,BCDX,BCEX, CDEX)}\mathbf{T}(X;Cq).
\end{align*} It gives the complete proof of \eqref{MMMMM}.

Obviously, when $X=q$, \eqref{MMMMM} reduces to
\begin{align}
\mathbf{T}(q;C)=\frac{\nabla\left(BCDEq^3,BC/q,CD/q,CE/q\right)}
{\nabla(C/q^3,BCDq,BCEq, CDEq)}\mathbf{T}(q;Cq).\label{watson-rec}
\end{align}
By iterating \eqref{watson-rec} $m$ times, we obtain
\begin{align}
\mathbf{T}(q;C)&=\frac{\poq{BCDEq^3,BC/q,CD/q,CE/q}{m}}
{\poq{C/q^3,BCDq,BCEq, CDEq}{m}}\mathbf{T}(q;Cq^m).\label{bbbbb}
\end{align}
Since $T(q;C)$ is analytic at $C=0$ and
\begin{align*}
\lim_{m\to +\infty}\mathbf{T}(q;Cq^m)=\mathbf{T}(q;0)=1,
\end{align*}
\eqref {bbbbb} reduces to
\begin{align*}
\mathbf{T}(q;C)=\frac{\poq{BCDEq^3,BC/q,CD/q,CE/q}{\infty}}
{\poq{C/q^3,BCDq,BCEq, CDEq}{\infty}}.
\end{align*} Thus we have \eqref{watson}. The lemma is proved.
\qed

We remark here that \eqref{watson} is just Rogers' ${}_6\phi_5$ summation formula  \citu{10}{(II. 20)}.  Even more, from Lemma \ref{MMMMM-YL} we obtain the following result.
\begin{yl}\label{MMMMM-YL-YL}Let $\mathbf{T}(X;C)$ be the same as in Lemma \ref{MMMMM-YL}. Then there exists certain function
$Q(X;B,D,E)$ being independent of $C$, such that
\begin{align}
\mathbf{T}(X;C)=Q(X;B,D,E)\frac{\poq{BCDEX^2q,q/(BCDEX^2),BC/q,CD/q,CE/q}{\infty}}
{\poq{1/(BCDEX), C/q^3,BCDX,BCEX, CDEX}{\infty}}.\label{NNNNN}
\end{align}
\end{yl}
\pf It only needs to consider the function
\begin{align}
F(C):=\frac{\poq{BCDEX^2q,q/(BCDEX^2),BC/q,CD/q,CE/q}{\infty}}
{\poq{1/(BCDEX), C/q^3,BCDX,BCEX, CDEX}{\infty}}
\end{align}
and  to check that
\[F(C)=\frac{\nabla\left(1/(BCDEXq),BCDEX^2q,BC/q,CD/q,CE/q\right)}
{\nabla(1/(BCDEX^2),C/q^3,BCDX,BCEX, CDEX)}F(Cq).\]
A direct comparison with \eqref{MMMMM} of Lemma  \ref{MMMMM-YL} yields
\begin{align}
\frac{\mathbf{T}(X;C)}{F(C)}=\frac{\mathbf{T}(X;Cq)}{F(Cq)}=\cdots=\frac{\mathbf{T}(X;Cq^m)}{F(Cq^m)}, m\geq 0
\end{align}
Next, we appeal to the uniqueness of Laurent series expansion, which states that if certain function $G(x)$ satisfies
\begin{align*}
G(x):=\sum_{n=-\infty}^\infty a_nx^n=\sum_{n=-\infty}^\infty a_n(xq)^n=G(xq),
\end{align*}
then there must hold that $a_n=0$ for $n\neq 0$. As such,  we now define such constant $a_0$ by
\begin{align}
Q(X;B,D,E):=\frac{\mathbf{T}(X;C)}{F(C)}.
\end{align}
Apparently, it is independent of $C$.
The lemma is thereby proved.
\qed

 We are now prepared  to show an equivalent variant  of Theorem \ref{baileypsi66}.
\begin{dl}[Bailey's  VWP ${}_6\psi_6$ summation formula]Let $B,C,D,E,X$  be
 arbitrary  nonzero complex parameters subject to $|C/q^3|<1$. Then there holds
\begin{align}
&{}_6\psi_6\bigg(BCDEX^2;BCDEXq,BXq,DXq,EXq;q,\frac{C}{q^3}\bigg)\label{XXXXX}
\\
&=\frac{\poq{q, 1/Bq,1/Dq, 1/Eq,BCDEX^2q,q/(BCDEX^2),BC/q,CD/q,CE/q}{\infty}}
{\poq{X, 1/BX,1/DX, 1/EX,1/(BCDEX), C/q^3,BCDX,BCEX, CDEX}{\infty}}.\nonumber
\end{align}
\end{dl}
\pf Obviously, by Lemma \ref{MMMMM-YL-YL} we  need only to find $Q(X;B,D,E)$. On account of its being independent of $C$, we now set in \eqref{NNNNN}  $BCDEXq=1$, namely,
$$C=\frac{1}{BDEXq}.$$ As is expected, we find
\begin{align*}
&Q(X;B,D,E)\\
&=\bigg(\mathbf{T}(X;C)\frac{\poq{1/(BCDEX), C/q^3,BCDX,BCEX, CDEX}{\infty}}{\poq{BCDEX^2q,q/(BCDEX^2),BC/q,CD/q,CE/q}{\infty}}\bigg)\bigg|_{C=\frac{1}{BDEXq}}.
\end{align*}
In this case, it is easy to check by the definitions of  ${}_r\phi_{r-1} $ and ${}_r\psi_r$ series (see \eqref{phino} and \eqref{yes}) that
\begin{align}
\mathbf{T}\bigg(X;\frac{1}{BDEXq}\bigg)
&={}_6\psi_6\bigg(X/q;1,BXq,DXq,EXq;q,\frac{1}{BDEXq^4}\bigg)\nonumber\\
&={}_6\phi_5\bigg(q/X;Bq^2,Dq^2,Eq^2;q,\frac{1}{BDEXq^4}\bigg).\label{watson-now0}
\end{align}
As already  proved, \eqref{watson}  asserts that for any $B,C,D,E$, it holds
\begin{align}
{}_6\phi_5(BCDEq^2;Bq^2,  Dq^2, Eq^2;q,C/q^3)=\frac{\poq{BCDEq^3,BC/q,CD/q,CE/q}{\infty}}
{\poq{C/q^3,BCDq,BCEq, CDEq}{\infty}}.\label{watson-now}
\end{align}
Therefore, by replacing
$C$ in \eqref{watson-now} with $1/(BDEXq)$ and substituting the result back to \eqref{watson-now0}, we obtain
\begin{align*}
\mathbf{T}\bigg(X;\frac{1}{BDEXq}\bigg)=\frac{\poq{q^2/X,1/(BDXq^2),1/(BEXq^2),1/(DEXq^2)}{\infty}}
{\poq{1/(BDEXq^4),1/BX,1/DX, 1/EX}{\infty}}.
\end{align*}
Finally, we arrive at
\begin{align*}
Q(X;B,D,E)&=\frac{\poq{q^2/X,1/(BDXq^2),1/(BEXq^2),1/(DEXq^2)}{\infty}}
{\poq{1/(BDEXq^4),1/BX,1/DX, 1/EX}{\infty}}\\
&\times \frac{\poq{q, 1/(BDEXq^4),1/Bq,1/Dq, 1/Eq}{\infty}}{\poq{X,q^2/X,1/(BDXq^2),1/(BEXq^2),1/(DEXq^2)}{\infty}}\\
&=\frac{\poq{q, 1/Bq,1/Dq, 1/Eq}{\infty}}{\poq{X, 1/BX,1/DX, 1/EX}{\infty}}.
\end{align*}
Upon substituting this back to \eqref{NNNNN}, we obtain \eqref{XXXXX} immediately. This completes our proof.
\qed

We conclude out paper with the following comments.
\begin{remark} It is easy to verify that \eqref{XXXXX}  reduces to  \eqref{bailey-a} at once by making the replacement
\begin{align*}
(B,C,D,E,X)\to\left(\frac{b c}{a q^2},\frac{a^2 q^4}{b c d e},\frac{b d}{a q^2},\frac{b e}{a q^2},\frac{a q}{b}\right).
\end{align*}
\end{remark}
\begin{remark} We may view the recurrence relation \eqref{bailey-rec} as saying that it is a common source and finite version for both Rogers'  ${}_6\phi_5$  and  Bailey's VWP ${}_6\psi_6$ summation formula.
\end{remark}

\begin{remark} It is clear that \eqref{trivalweierstrass-new-1} of Lemma \ref{theta} is in fact the special case $q=0$ of the famous Weierstrass  theta identity (see \citu{10}{Exercise 2.16(i)} or \cito{koornwinder}{})
\begin{align}
\theta\left(cx,\frac{x}{c},bz,
\frac{z}{b};q\right)-\theta\left(bx,\frac{x}{b},cz,\frac{z}{c};q\right)=\frac{z}{c}\theta\left(bc,\frac{c}{b},
xz,\frac{x}{z};q\right).\label{thetaidd-newnew}
\end{align}
Note that $\theta(x;q)$ denotes the  Jacobi modified theta function given by
\begin{align*}
\poq{x,q/x}{\infty}
\end{align*}
and its multi-parameter form
\[\theta(x_1,x_2,\cdots,x_m;q):\:=\:\theta(x_1;q)\theta(x_2;q)\cdots\theta(x_m;q).\]
The the reader may consult \cito{koornwinder}{} for  a full history and further applications concerning
Weierstrass' theta identity. It is  worth mentioning  that in the recent paper \cite[Theorem 1.7]{kn:wangjinpaper}, the author showed that \eqref{thetaidd-newnew}   is equivalent to  \eqref{trivalweierstrass-new-1}.
\end{remark}

\end{document}